\documentclass{article}

\usepackage{amsmath, amssymb, amsfonts, amsthm}
\usepackage{graphicx}
\usepackage{color}
\usepackage[colorlinks,linkcolor=blue,citecolor=red,urlcolor=blue]{hyperref}
\usepackage{tikz}

\numberwithin{equation}{section}

\theoremstyle{plain}
\newtheorem{thm}{Theorem}[section]

\theoremstyle{definition}

\newtheorem{rem}[thm]{Remark}

\begin{document}

% ----------------------------------------------------------------------
% Title and Author Information
% ----------------------------------------------------------------------
\title{\textbf{Constraint Structure and Zero Counting in the Integral Representation of the Zeta Function}}

\author{
    \textbf{Nianrong Feng} \\
    \small{School of Computer and Information, Anhui Normal University}\\
%    \small{Wuhu 241000, China}\\
    \small{\texttt{fengnr@mail.ahnu.edu.cn}}
    \and
    \textbf{Yongzheng Wang} \\
    \small{School of Mathematics and Statistics, Anhui Normal University}\\
%    \small{Wuhu 241000, China}\\
    \small{\texttt{nrfeng@mail.ahnu.edu.cn}}
}

\date{}

\maketitle

% ----------------------------------------------------------------------
% Abstract (Updated with "Two-End Anchoring", "Construct", "Grounded in")
% ----------------------------------------------------------------------
\begin{abstract}
\noindent Starting from the classical integral representation of the $\zeta(s)$ function introduced by Riemann in 1859, this paper reexamines its analytic symmetry structure. By performing a geometric decomposition of the integral representation, we demonstrate that on the critical line $\Re(s)=\frac{1}{2}$, the value of $\xi(s)$ corresponds strictly to the \textbf{real-part projection} of a specific analytic component. This discovery equivalently transforms the problem of complex zeros into a problem of \textbf{sign evolution} along the real axis.

Based on this geometric framework, we \textbf{construct} an analytic mechanism of \textbf{"Two-End Anchoring, Interval Counting"}: the global argument increment on the region boundary \textbf{anchors} the initial value of the phase function, while the geometric decomposition structure on the critical line \textbf{locks} its final value. This mechanism reveals an \textbf{intrinsic coherence} between global topological constraints and local sign oscillations. Unlike traditional methods that rely on asymptotic estimates (such as the Big $O$ error term), the analysis in this paper is \textbf{grounded in} exact identities. It unveils the \textbf{geometric determinism} underlying the zero-counting formula, offering a novel perspective for analytic number theory independent of asymptotic analysis.

\vspace{0.5cm}

\end{abstract}

% ----------------------------------------------------------------------
% Main Content
% ----------------------------------------------------------------------

\section{Introduction}

In 1859, in his seminal paper \textit{"On the Number of Primes Less Than a Given Magnitude"} \cite{riemann}, Riemann systematically introduced complex variable methods to study the distribution of prime numbers for the first time. Historically, this work did not emerge in isolation within analytic number theory but was deeply rooted in Riemann's holistic research on complex analysis and geometric theory at the time. In the years preceding this, Riemann had systematically developed a methodology in complex analysis (Riemann surface theory) \cite{riemann1851}, geometry (his famous 1854 inaugural lecture) \cite{riemann1854}, as well as in the theory of elliptic functions \cite{riemann1857} and Abelian functions \cite{riemann1857}, treating analytic objects as geometric entities \cite{riemann1857, laugwitz}. Within this conceptual framework, introducing the $\zeta(s)$ function and examining its zero structure can be viewed as a natural extension of his geometric-analytic thought into number theory \cite{laugwitz}.

In Riemann's approach, $\zeta(s)$ was not studied as an isolated analytic function but was placed within an analytic framework possessing explicit symmetry through its integral representation and functional equation. The functional equation guarantees the analytic continuation of $\zeta(s)$ on the one hand, and establishes the critical line $\Re(s)=1/2$ as a fundamental axis of symmetry in the geometric sense on the other. It is particularly noteworthy that the integral representation introduced by Riemann in his paper:
\begin{equation}\label{sggsx1}
\pi^{-\frac{s}{2}}\Gamma\left(\frac{s}{2}\right)\zeta(s)
= \frac{1}{s(s-1)} + \int_1^{\infty}\psi(x)
\left(x^{\frac{s}{2}-1}+x^{-\frac{s+1}{2}}\right)\,dx,
\end{equation}
was not merely used for deriving the functional equation. Structurally, this representation already contains a symmetric decomposition: when $s$ takes the form of the critical line $s=\tfrac12+it$, the two integral terms on the right-hand side are complex conjugates of each other, thereby naturally decomposing the entire expression into the sum of an analytic function and its conjugate. This property is not an additional assumption but a result generated by the combined action of the kernel's symmetry and the functional equation's structure, playing a decisive role in the subsequent analysis of the behavior on the critical line.

In the ensuing century and a half, research on the zero distribution of the $\zeta(s)$ function gradually followed the main routes of asymptotic counting and numerical computation. Specifically, the Riemann-von Mangoldt formula \cite{Titch} revealed the global counting law of zeros through argument variation; meanwhile, research focusing on the critical line $\Re(s)=\frac{1}{2}$ centered on local oscillatory behavior, zero density, and numerical verification, with the Hardy-Littlewood method \cite{Titch} and the Riemann-Siegel formula \cite{xref1} being typical representatives. Although these methods achieved significant success at the computational level, they also weakened the geometric information inherent in the functional equation and its original integral representation.

Of particular note is the identity explicitly recorded by Siegel when organizing Riemann's Nachlass (commonly referred to as Formula (60)) \cite{xref1}. For a long time, it primarily served as a starting point for numerical approximation and saddle-point analysis, while its precise geometric significance was rarely explored in depth. In the traditional research framework, the global zero-counting mechanism and the argument structure on the critical line were often analyzed separately, and the intrinsic correlation between the two has not yet been systematically revealed. This separation largely stems from the reliance on asymptotic counting formulas, which in turn obscures the global geometric constraints implied by exact identities.

This study returns to the original form of Riemann's integral representation. Without relying on asymptotic estimates, we examine the geometric argument constraints implied by the functional equation on the critical line $\Re(s)=\frac{1}{2}$. The focus of this research is not on obtaining new numerical results, but on clarifying the structural correlations inherent between different zero-counting methods. From this perspective, the zero problem is no longer understood merely as an analytic condition for a complex function to vanish, but as the intersection positions that naturally emerge during the argument evolution under global geometric structural constraints.

It is worth pointing out that the research perspective adopted in this paper embodies a profound shift from "asymptotic estimation" to "structural determination." This shift can be illustrated by an analogy from elementary analysis: the fact that $\sin(\frac{\pi}{6})=\frac{1}{2}$ is essentially an intrinsic consequence of the geometric symmetry of the unit circle, rather than merely the asymptotic result of the Taylor series $x-\frac{x^3}{3!}+\dots$ at a specific point. Similarly, this paper aims to demonstrate that the \textbf{zero distribution of the $\zeta(s)$ function on the critical line} is not just a statistical phenomenon in the sense of asymptotic analysis, but a structural result determined by the \textbf{Geometric Rigidity} jointly established by the functional equation and the integral representation.

In terms of methodology, this paper does not directly track the complex oscillatory behavior of $\zeta(s)$ on the critical line via analytic means. Instead, it establishes a geometric mechanism of \textbf{"Two-End Anchoring, Interval Counting."} Specifically, we introduce the geometric \textbf{phase function} $\omega(t)$ naturally induced by the integral representation and transform the zero-counting problem into a one-dimensional geometric counting process: the global zero-counting formula within the region \textbf{anchors} the initial value of $\omega(t)$, while the special status of the point $s=\tfrac12+iT$ in the geometric decomposition \textbf{locks} its final value. Within this rigid framework, the number of zeros on the critical line is strictly equivalent to the number of times the \textbf{phase function} $\omega(t)$ crosses odd multiples of $\tfrac{\pi}{2}$ within this determined range. This transformation elevates the zero problem from an analytic task dependent on local error control to a geometric counting problem decided jointly by the global symmetric structure and the continuity principle, thereby restoring the original geometric nature of Riemann's thought in number theory without relying on asymptotic expansions.

\section{Notations and Basic Objects}

Let $\zeta(s)$ denote the analytic continuation of the Riemann zeta function. Define
\[
\theta(t)=\arg\!\left[\pi^{-it/2}\Gamma\!\left(\frac14+\frac{it}{2}\right)\right].
\]
For $T>0$, we define the rectangular region $\mathbf D(T)$ and the critical line segment $\mathbf{L}(T)$ respectively as:
\[
\mathbf D(T)=\{s=\sigma+it \in \mathbb{C} :-1\le\sigma\le 2,\; 0\le {t}\le T\},
\]
\[
\mathbf{L}(T) = \{s=\sigma+it \in \mathbb{C} : \sigma= 1/2, \, 0\leq t \leq T\}.
\]
For convenience of discussion, without loss of generality, we assume that there are no zeros of $\zeta(s)$ on the boundary of $\mathbf D(T)$.

\section{Geometric Framework and Structural Symmetry}

\subsection{Global Symmetry Induced by the Functional Equation}

In the analytic study of the $\zeta(s)$ function, the symmetry embodied by the functional equation plays a fundamental role. To this end, we introduce the following standard notations for entire functions:
\[
\Phi(s)=\pi^{-s/2}\Gamma\!\left(\frac{s}{2}\right)\zeta(s),
\qquad
\xi(s)=\tfrac12 s(s-1)\Phi(s).
\]
The classical Riemann functional equation
\begin{equation}\label{eq:functional_eq}
\xi(s)=\xi(1-s)
\end{equation}
indicates that $\xi(s)$ possesses strict reflection symmetry with respect to the line $\Re(s)=\tfrac12$. Consequently, if $\xi(s)$ has zeros in the critical strip that do not lie on the critical line, these zeros must appear in pairs symmetric with respect to the critical line.

This symmetric structure makes the Argument Principle a natural tool for studying zero distribution. By examining the argument variation of $\xi(s)$ along the boundary of the rectangular region $\mathbf{D}(T)$, we can obtain global topological information about the zero count within that region.

\subsection{Geometric Decomposition and Real-Part Projection Mechanism}

The integral representation \eqref{sggsx1} introduced by Riemann in his 1859 paper not only realizes the analytic continuation of $\zeta(s)$ but also implies an \textbf{exact geometric decomposition} of the functional equation's structure. Specifically, there exists the following decomposition form:
\begin{equation}\label{implicit_sym}
\pi^{-\frac{s}{2}}\Gamma\!\left(\frac{s}{2}\right)\zeta(s)
=\varphi(s)+\varphi(1-s),
\end{equation}
where
\[
\varphi(s)=-\frac{1}{s}
+\int_1^{\infty}\psi(x)x^{\frac{s}{2}-1}\,dx .
\]
This decomposition was explicitly extracted by Siegel in his systematic organization of Riemann's Nachlass and appears in the form of the following exact identity (corresponding to Formula (60) in Siegel's work):
\begin{equation}\label{eq:siegel60}
\pi^{-s/2}\Gamma\!\left(\frac{s}{2}\right)\zeta(s)
=\varphi(s)+\varphi(1-s).
\end{equation}
In this representation, $\varphi(s)$ can be defined by a contour integral with a specific orientation:
\begin{equation}\label{eq:phi-def}
\varphi(s)
=\pi^{-s/2}\Gamma\!\left(\frac{s}{2}\right)
\int_{0\swarrow 1}
\frac{x^{-s}e^{\pi i x^2}}{e^{\pi i x}-e^{-\pi i x}}\,dx .
\end{equation}

It must be emphasized that the identities \eqref{implicit_sym} and \eqref{eq:siegel60} hold strictly throughout the entire complex plane. They directly reflect the intrinsic symmetric structure of the functional equation, rather than any asymptotic expansion or approximate expression.

When the variable $s$ is restricted to the critical line $s=\tfrac12+it$, using the conjugate symmetry of the kernel, it is easy to prove:
\[
\varphi(1-s)=\overline{\varphi(s)}.
\]
Therefore, the above geometric decomposition degenerates into the following \textbf{projection formula} on the critical line:
\begin{equation}\label{eq:real-projection}
\Phi\!\left(\tfrac12+it\right)
=2\,\Re\!\left[\varphi\!\left(\tfrac12+it\right)\right].
\end{equation}
This formula reveals a core geometric fact: on the critical line, the value of $\Phi(s)$ is completely determined by the real-part trajectory of the holomorphic function $\varphi(s)$. In other words, the oscillation of $\Phi(\tfrac12+it)$ is not an independent random behavior in the complex plane but is strictly controlled by the argument evolution of $\varphi(\tfrac12+it)$ through \textbf{orthogonal projection onto the real axis}.

\subsection{Argument Dynamics and Geometric Characterization of Zeros}

Based on the real-part projection mechanism described above, we define the \textbf{Phase Evolution Function} (or simply \textbf{Phase Function}) on the critical line:
\[
\omega(t)=\arg \varphi\!\left(\tfrac12+it\right),
\]
where the argument $\arg$ takes a continuous branch along the critical line. From the projection formula \eqref{eq:real-projection}, the zero problem of $\zeta(\tfrac12+it)$ on the critical line can be completely reduced to the threshold crossing problem of the phase function $\omega(t)$. Specifically, the real part $\Phi(\tfrac12+it)$ undergoes a sign flip (or tangency) if and only if $\omega(t)$ satisfies the condition:
\[
\omega(t)=\left(k+\tfrac12\right)\pi,\qquad k\in\mathbb Z,
\]
which corresponds to the appearance of a zero. Therefore, examining the zero distribution along the critical line is geometrically equivalent to studying the number of times the continuous function $\omega(t)$ crosses half-integer multiples of $\pi$ within a given interval.

This characterization avoids complex estimations of the modulus and asymptotic terms of $\zeta(\tfrac12+it)$, successfully transforming the complex zero problem into a \textbf{geometric counting problem} regarding the continuous real function $\omega(t)$. This one-to-one correspondence between argument evolution and real-part sign change constitutes the basis for the subsequent theoretical analysis in this paper.

\subsection{Intrinsic Coherence: From Local Argument to Global Topology}

The classical Argument Principle provides a zero-counting method that relies solely on the argument variation along the boundary $\partial\mathbf{D}(T)$, and its conclusion possesses topological stability. However, traditional views often regard the local behavior on the critical line as a statistical process independent of the boundary counting.

The geometric framework established by formula \eqref{eq:siegel60} indicates that when the discussion is restricted to the critical line $\Re(s)=\tfrac12$, the argument evolution of the function is strictly constrained by the global structure. Since $\Phi(\tfrac12+it)$ is essentially the projection of the analytic component $\varphi$, the rhythm of its sign changes (i.e., the generation of zeros) is directly governed by the phase evolution rate of $\varphi$.

Thus, there exists an \textbf{Intrinsic Coherence} between the local argument variation pattern on the critical line and the global argument increment on the region boundary. These two are not independent counting processes but are projection manifestations of the same analytic structure (i.e., the geometric rigidity induced by the functional equation) in different dimensions. The following two theorems will precisely characterize and prove this profound constraint relation.

\section{Core Geometric Theorems}

The total variation of the argument of the function $\xi(s)$ along the region boundary $\partial\mathbf{D}(T)$ is dictated solely by the boundary value $\arg \zeta(\tfrac12+iT)$ and the argument function $\theta(T)$ contributed by the $\Gamma(s)$ factor.

\begin{thm}[Global Geometric Constraint Theorem]\label{xdl2}
Inside the rectangular region $\mathbf{D}(T)$, the total number of zeros of $\xi(s)$, denoted by $N(T)$, is given by the following exact argument expression:
\begin{align}\label{NTF1}
N(T) = \frac{\theta(T) + \arg \zeta(\frac{1}{2}+iT)}{\pi} + 1,
\end{align}
where $\theta(T) = \arg[\pi^{-iT/2}\Gamma(\frac{1}{4}+\frac{iT}{2})]$.

\textbf{Remark.} This formula is classically expressed in the asymptotic form known as the Riemann-von Mangoldt formula:
\[
N(T)=\frac{T}{2\pi}\log\frac{T}{2\pi}-\frac{T}{2\pi}+O(\log{T}).
\]
However, to establish a rigorous topological constraint, it is imperative to distinguish the \textbf{exact argument expression} employed here from the traditional asymptotic approximation.
\end{thm}

As $t$ varies continuously, the evolution of the phase function $\omega(t)$ is strictly constrained by the geometric structure. Specifically, the sign changes of $\Re\varphi(\tfrac12+it)$ are restricted to occur only when $\omega(t)$ crosses odd multiples of $\frac{\pi}{2}$, precluding any \textbf{arbitrary local fluctuations}.

\begin{thm}[Phase Locking Theorem]\label{xdl3}
Based on the geometric structure jointly induced by the integral representation \eqref{sggsx1} and the functional equation, the total variation of the phase $\omega(t)$ along the critical line segment $\mathbf{L}(T)$ is uniquely determined by $\arg \zeta(\tfrac12+iT)$ and the argument function $\theta(T)$ derived from the $\Gamma(s)$ factor.
\end{thm}

\section{Proofs of Main Theorems}

\subsection{Proof of Theorem \ref{xdl2}}

\begin{proof}
Let $R$ denote the positively oriented rectangular contour of $D(T)$. Since the entire function $\xi(s)$ has no zeros on $R$, by the Argument Principle, the number of zeros of $\xi(s)$ inside $R$ is
\begin{align}\label{hljfa}
N(T)=\frac{1}{2\pi}\triangle_R\arg{\xi(s)},
\end{align}
where $\triangle_R\arg{\xi(s)}$ denotes the change in the argument of $\xi(s)$ along the contour $R$. Divide the contour $R$ into three segments: let $L_1$ be the horizontal line from $-1$ to $2$; $L_2$ be the polygonal line from $2$ to $2+iT$ and then to $\frac{1}{2}+iT$; and $L_3$ be the polygonal line from $\frac{1}{2}+iT$ to $-1+iT$ and then to $-1$. Therefore,
\[
\triangle_R\arg{\xi(s)}=\triangle_{L_1}\arg{\xi(s)}+\triangle_{L_2}\arg{\xi(s)}+\triangle_{L_3}\arg{\xi(s)}.
\]
From $\xi(s)=\xi(1-s)$ and $\xi(s)=\overline{\xi(1-\overline{s})}$, we obtain
\[
\triangle_{L_2}\arg{\xi(s)}=\triangle_{L_3}\arg{\xi(s)}.
\]
Since $\xi(s)$ is real on $L_1$, $\triangle_{L_1}\arg{\xi(s)}=0$. Thus,
\begin{align*}
N(T)&=\frac{1}{\pi}\triangle_{L_2}\arg{\xi(s)}.
\end{align*}
We have
\begin{align*}
\triangle_{L_2}\arg{\xi(s)}&=\triangle_{L_2}\arg{\left[\frac{s}{2}(s-1)\pi^{-s/2}\Gamma(\frac{s}{2})\zeta(s)\right]}\\
&=\triangle_{L_2}\arg{(\frac{s}{2})}+\triangle_{L_2}\arg{(s-1)}+\triangle_{L_2}\arg{\left[\pi^{-s/2}\Gamma(\frac{s}{2})\zeta(s)\right]}.
\end{align*}
Considering $\triangle_{L_2}\arg{(\frac{s}{2})}$ and $\triangle_{L_2}\arg{(s-1)}$ separately:
\[
\triangle_{L_2}\arg{(\frac{s}{2})}=\arg{\left(\frac{\frac{1}{2}+iT}{2}\right)}-\arg{1}=\arg{\left(\frac{1}{2}+iT\right)},
\]
\[
\triangle_{L_2}\arg{(s-1)}=\arg{\left(-\frac{1}{2}+iT\right)}-\arg{1}=\arg{\left(-\frac{1}{2}+iT\right)}.
\]
The argument change of $\frac{s}{2}(s-1)$ on $L_2$ is
\begin{align*}
\triangle_{L_2}\arg{\left[\frac{s}{2}(s-1)\right]}&=\triangle_{L_2}\arg{(\frac{s}{2})}+\triangle_{L_2}\arg{(s-1)}\\
&=\arg{\left(\frac{1}{2}+iT\right)}+\arg{\left(-\frac{1}{2}+iT\right)}\\
&=\pi.
\end{align*}
When $s=2$, $\pi^{-s/2}\Gamma(\frac{s}{2})\zeta(s)$ is a real function, so its argument is 0. Thus,
\begin{align*}
\triangle_{L_2}\arg{\left[\pi^{-s/2}\Gamma(\frac{s}{2})\zeta(s)\right]}&=\arg{\left[\pi^{-\frac{1/2+iT}{2}}\Gamma(\frac{1/2+iT}{2})\zeta(1/2+iT)\right]}-0\\
&=\arg{\left[\pi^{-\frac{1/2+iT}{2}}\Gamma(\frac{1/2+iT}{2})\right]}+\arg{\zeta(1/2+iT)}\\
&=\theta(T)+\arg{\zeta(1/2+iT)}.
\end{align*}
Therefore,
\begin{align*}
N(T)&=\frac{1}{\pi}\triangle_{L_2}\arg{\xi(s)}\\
&=\frac{1}{\pi}\triangle_{L_2}\arg{\left[\pi^{-s/2}\Gamma(\frac{s}{2})\zeta(s)\right]}+\frac{1}{\pi}\triangle_{L_2}\arg{\left[\frac{s}{2}(s-1)\right]}\\
&=\frac{\theta(T)+\arg\zeta(1/2+iT)}{\pi}+1.
\end{align*}
This completes the proof.
\end{proof}

\subsection{Proof of Theorem \ref{xdl3}}

\begin{proof}
When $s=\frac{1}{2}+it$, from equation \eqref{implicit_sym} we obtain
\begin{align}\label{llld}
{\pi^{-\frac{s}{2}}\Gamma(\frac{s}{2})\zeta(s)}&=2\Re{\left[\varphi(s)\right]}\nonumber\\
&={2r(t)\cos[\omega(t)]},
\end{align}
where $r(t)=|\varphi(s)|$ and $\omega(t)=\arg\varphi(s)$.

\begin{figure}[htbp]
  \centering
  % Note: Ensure the image file 'fourgg.png' is in the same directory.
  \includegraphics[width=0.8\textwidth]{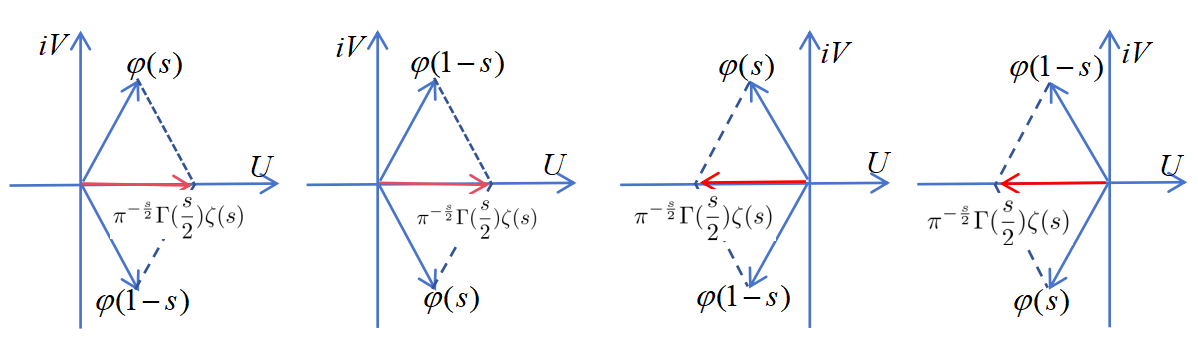}
  \caption{Four possible geometric configurations of equation \eqref{implicit_sym} on the critical line $s=1/2+iT$, demonstrating the real-part projection mechanism induced by the functional equation.}
  \label{fig:fourgg}
\end{figure}

When the variable takes values on the critical line $s=1/2+iT\in{\mathbf{D}(T)}$, the geometric configuration induced by equation \eqref{implicit_sym} is shown in Figure \ref{fig:fourgg}. According to this real-part projection mechanism, the following argument relation between the auxiliary component $\varphi(s)$ and the composite vector $\pi^{-\frac{s}{2}}\Gamma(\frac{s}{2})\zeta(s)$ can be directly derived:
\begin{align*}
\left|\arg\varphi(s)- \arg\pi^{-\frac{s}{2}}\Gamma(\frac{s}{2})\zeta(s)\right|<{\frac{\pi}{2}}.
\end{align*}
Assuming non-zero magnitude, we have:
\begin{align}\label{fsx}
-\frac{\pi}{2}<\arg\varphi(s)-\arg\pi^{-\frac{s}{2}}\Gamma(\frac{s}{2})-\arg\zeta(s)<{\frac{\pi}{2}}.
\end{align}
According to equation \eqref{llld}, the function $\zeta(1/2+it)$ has two types of zeros:
\begin{enumerate}
    \item $N_{\cos}(T)$: The number of zeros satisfying $\cos[\omega(t)]=0$.
    \item $N_{r}(T)$: The number of zeros satisfying $r(t)=0$ and $\cos[\omega(t)]\neq 0$.
\end{enumerate}
Thus, the total number of zeros of $\zeta(1/2+it)$ on $L(T)$ is
\begin{align}\label{N0Tsum}
N_0(T)=N_{\cos}(T)+N_r(T).
\end{align}

We examine the zero distribution of the function $\cos[\omega(t)]$ on the interval $t\in[0,T]$. As $t$ varies, the phase function $\omega(t)$ evolves continuously from $\omega(0)$ to $\omega(T)$. Without making additional assumptions about the specific evolution of $\omega(t)$, we can assert that when $\omega(t)$ is strictly monotonically increasing on $[0,T]$, the number of zeros of $\cos[\omega(t)]$ attains its minimum value; any local non-monotonic oscillation would result in an increase in the number of zeros. Therefore, the following analysis is primarily based on the monotonic case to obtain a baseline estimate for the zero count.

\medskip
\textbf{Step 1: Determining the initial argument interval value}

It is known that $t_1\doteq14.13$ is the first nontrivial zero of $\zeta(s)$. According to Theorem \ref{xdl2}, when $T<t_1$, the rectangular region $\mathbf D(T)$ contains no zeros, i.e., $N(T)=0$. Taking $T=0(<t_1)$ and substituting it into the zero-counting formula \eqref{NTF1}, we obtain
\[
N(0)=\frac{\theta(0)+\arg \zeta(\tfrac12)}{\pi}+1=0,
\]
which implies
\begin{equation}\label{tmpf1}
\theta(0)+\arg \zeta(\tfrac12)=-\pi .
\end{equation}

On the other hand, substituting $T=0\in\mathbf D(T)$ into inequality \eqref{fsx} yields
\[
-\frac{\pi}{2}
<
\arg\varphi(\tfrac12)
-
\theta(0)
-
\arg \zeta(\tfrac12)
<
\frac{\pi}{2}.
\]
Combining this with equation \eqref{tmpf1}, the above expression simplifies to
\[
-\frac{3\pi}{2}
<
\arg\varphi(\tfrac12)
<
-\frac{\pi}{2}.
\]
Substituting $\omega(t)=\arg\varphi\!\left(\tfrac12+it\right)$ into the inequality, we obtain the value interval for the initial phase:
\begin{equation}\label{omegabds1}
-\frac{3\pi}{2}<\omega(0)<-\frac{\pi}{2}.
\end{equation}

\medskip
\textbf{Step 2: End argument interval and baseline zero count}

It is known that $t_2\doteq21.02$ and $t_3\doteq25.01$ are the second and third nontrivial zeros of $\zeta(s)$, respectively. Take any
\[
T_1\in(t_2,t_3).
\]
From Theorem \ref{xdl2}, we have $N(T_1)=2$. Substituting $T_1$ into \eqref{NTF1}, we get
\[
\frac{\theta(T_1)+\arg \zeta(\tfrac12+iT_1)}{\pi}+1=2,
\]
which simplifies to
\begin{equation}\label{tmpf2}
\theta(T_1)+\arg \zeta(\tfrac12+iT_1)=\pi .
\end{equation}
Substituting $T_1$ into inequality \eqref{fsx} gives
\[
-\frac{\pi}{2}
<
\arg\varphi(\tfrac12+iT_1)
-
\theta(T_1)
-
\arg \zeta(\tfrac12+iT_1)
<
\frac{\pi}{2}.
\]
Combining this with equation \eqref{tmpf2}, we obtain
\[
\frac{\pi}{2}
<
\arg\varphi(\tfrac12+iT_1)
<
\frac{3\pi}{2},
\]
that is,
\begin{equation}\label{omegabds2}
\frac{\pi}{2}<\omega(T_1)<\frac{3\pi}{2}.
\end{equation}
From \eqref{omegabds1} and \eqref{omegabds2}, using the continuity of $\omega(t)$ with respect to $t$, by the Intermediate Value Theorem, there must exist $T_0\in(0,T_1)$ such that
\[
\omega(T_0)=0.
\]
Using $T_0$ as a dividing point, we decompose the interval $[0,T]$ into
\[
[0,T]=[0,T_0)\cup[T_0,T].
\]
Correspondingly, the phase interval satisfies
\[
[\omega(0),\omega(T)]
=
[\omega(0),0)\cup[0,\omega(T)].
\]

\medskip
\noindent
\textbf{(i) Zeros in the interval $\boldsymbol{\omega(t)\in[\omega(0),0)}$}

Examining equation \eqref{omegabds1}, it is certain that $-\tfrac{\pi}{2}\in(\omega(0),0)$. Thus, there exists $t_c\in(0,T_0)$ such that
\[
\omega(t_c)=-\frac{\pi}{2},
\qquad
\cos\omega(t_c)=0.
\]
This indicates that within the interval $\omega(t)\in[\omega(0),0)$, the function $\cos[\omega(t)]$ has at least one zero, corresponding to a zero $s=\tfrac12+it_c$ on the critical line.

\medskip
\noindent
\textbf{(ii) Zeros in the interval $\boldsymbol{\omega(t)\in[0,\omega(T)]}$}

Note that
\[
\cos\!\left(k\pi-\tfrac{\pi}{2}\right)=0\qquad(k\in\mathbb Z).
\]
Therefore, in the interval $\omega(t)\in[0,\omega(T)]$, the number of zeros of $\cos[\omega(t)]$ is
\[
k=\left\lfloor\frac{\omega(T)}{\pi}+\frac{1}{2}\right\rfloor .
\]
To determine the value of $\omega(T)$, utilizing inequality \eqref{fsx} again, we have
\[
\frac{\theta(T)+\arg\zeta(\tfrac12+iT)}{\pi}
<
\frac{\omega(T)}{\pi}+\frac{1}{2}
<
\frac{\theta(T)+\arg\zeta(\tfrac12+iT)}{\pi}+1.
\]
Thus,
\[
\left\lfloor\frac{\omega(T)}{\pi}+\frac{1}{2}\right\rfloor
=
\frac{\theta(T)+\arg\zeta(\tfrac12+iT)}{\pi}.
\]
In summary, the total number of zeros of $\cos[\omega(t)]$ in the interval $[\omega(0),0)\cup[0,\omega(T)]$ is
\[
N_{\cos}(T)
=k+1=
\frac{\theta(T)+\arg\zeta(\tfrac12+iT)}{\pi}+1.
\]
Substituting $N_{\cos}(T)$ into formula \eqref{N0Tsum}, we obtain
\[
N_0(T)
=
\frac{\theta(T)+\arg\zeta(\tfrac12+iT)}{\pi}+1+N_r(T).
\]
Combining this with \eqref{NTF1}, we have
\[
N_0(T)=N(T)+N_r(T).
\]
Since it always holds that $N(T)\ge N_0(T)$ and $N_r(T)\ge0$, it must be that $N_r(T)=0$. Therefore, within this geometric structural framework, the zero count on the critical line remains consistent with the zero count in the region.
\end{proof}

\begin{rem}
It is crucial to note that the conclusions above do not rely on any \textit{a priori} assumptions regarding the monotonicity or local regularity of the phase function
\(
\omega(t)=\arg\varphi(\tfrac12+it)
\).
We have not artificially excluded potential local oscillations at the outset of our analysis.

The constraints imposed on the phase evolution along the critical line arise not from additional analytic hypotheses, but directly from the \textbf{structural rigidity} inherent in the exact identity \eqref{eq:siegel60} and its real-part projection mechanism \eqref{eq:real-projection}. This mechanism mandates that the local phase evolution must maintain strict \textbf{geometric coherence} with the global topological invariants determined by the boundary of $\mathbf D(T)$.

Consequently, any "free" local oscillation that deviates from this global constraint is automatically precluded at the structural level, as it would violate the logical consistency between the sign changes of $\Phi(\tfrac12+it)$ and the global conservation of argument.
\end{rem}

\section{Dependence Relation Between Zero Counts}

In the research framework of classical analytic number theory, global zero counting within a region and zero distribution on the critical line are usually studied as two independent technical subjects. The former relies on the Argument Principle, using the argument variation on the region boundary to obtain the total number of zeros; the latter focuses on the local asymptotic analytic behavior of the function near the critical line. Although this artificial distinction provides operational convenience for asymptotic analysis, it fundamentally severs the intrinsic structural connection between the two counting methods.

The geometric decomposition framework established based on Riemann's original integral representation \eqref{eq:siegel60} reveals the other side of the fact: the argument variation on the region boundary and the phase evolution on the critical line do not originate from mutually independent analytic mechanisms. When the variable is restricted to the critical line $\Re(s)=\tfrac12$, the value of $\Phi(s)$ is \textbf{forcibly locked} as the real-part projection of the auxiliary analytic component $\varphi(s)$, thereby making the zero counting on the critical line structurally strictly \textbf{subordinate} to the global topological information determined by the region boundary.

Therefore, the zero distribution on the critical line should not be understood as some independent random oscillation, but as the \textbf{geometric projection} of the global region counting on the symmetry axis of the functional equation. This deep dependence relation is established entirely on exact identities and the continuity principle; although it is often obscured by error terms in traditional asymptotic analysis, in the exact geometric framework of this paper, it manifests as an \textbf{unavoidable mathematical inevitability}.

\section{Historical Perspective and Methodological Reflection}

It is particularly worth reflecting that the identity \eqref{eq:siegel60}, explicitly extracted by Siegel when organizing Riemann's Nachlass, clearly presented the symmetric decomposition structure regarding $\varphi(s)+\varphi(1-s)$ in form. However, for a long time, this formula has primarily been regarded as the algorithmic basis for numerical computation (the Riemann-Siegel formula), while the \textbf{Geometric Rigidity} and \textbf{Structural Duality} inherent behind it have not received due theoretical attention. In a sense, the work of this paper is a return to Riemann's original geometric intuition.

\section{Conclusion}

From a methodological perspective, the analytical framework adopted in this paper can be summarized as a geometric counting mechanism dominated by the global symmetric structure. Through the exact identities established by the Riemann integral representation and the functional equation, the zero problem on the critical line is effectively transformed into an interval counting problem regarding the geometric phase function.

Specifically, the global zero-counting formula within the region provides a natural anchor for the initial value of the phase function, while the special status of the endpoint $s=\tfrac12+iT$ in the geometric decomposition determines its final value. Under this rigid framework constrained by both ends, the number of zeros on the critical line is completely determined by the crossing behavior of the phase function within this determined range, no longer relying on fine control of the function's local oscillations.

This perspective indicates that zero counting on the critical line and zero counting in the region are not two independent analytical means, but manifestations of the same analytic-geometric structure at different levels. The results presented in this paper do not rely on asymptotic expansions or error term estimations, but are directly established on exact identities and the continuity principle, thereby offering a new analytical perspective for understanding the structural consistency in the zero distribution of the Riemann $\zeta$ function.


\begin{thebibliography}{99}

\bibitem{riemann}
B. Riemann, \textit{\"{U}ber die Anzahl der Primzahlen unter einer gegebenen Gr\"{o}sse}, Monatsberichte der Berliner Akademie, 1859.

\bibitem{Siegel}
C. L. Siegel, \textit{\"{U}ber Riemanns Nachlass zur analytischen Zahlentheorie}, Quellen und Studien zur Geschichte der Mathematik, Astronomie und Physik, \textbf{2} (1932), 45--80. (Also in \textit{Gesammelte Abhandlungen}, Bd. I, Springer, 1966).

\bibitem{riemann1851}
B. Riemann, \textit{Grundlagen f\"{u}r eine allgemeine Theorie der Functionen einer ver\"{a}nderlichen complexen Gr\"{o}sse}, Dissertation, G\"{o}ttingen, 1851.

\bibitem{riemann1854}
B. Riemann, \textit{\"{U}ber die Hypothesen, welche der Geometrie zu Grunde liegen}, Abhandlungen der K\"{o}niglichen Gesellschaft der Wissenschaften zu G\"{o}ttingen, \textbf{13} (1867), 133--152. (Habilitationsvortrag, 1854).

\bibitem{riemann1857}
B. Riemann, \textit{Theorie der Abel'schen Functionen}, Journal f\"{u}r die reine und angewandte Mathematik, \textbf{54} (1857), 115--155.

\bibitem{laugwitz}
D. Laugwitz, \textit{Bernhard Riemann 1826--1866: Turning Points in the Conception of Mathematics}, Birkh\"{a}user, Basel, 1999.

\bibitem{Titch}
E. C. Titchmarsh, \textit{The Theory of the Riemann Zeta-Function}, 2nd ed. (Revised by D. R. Heath-Brown), Clarendon Press, Oxford, 1986.

\bibitem{xref1}
H. M. Edwards, \textit{Riemann's Zeta Function}, Academic Press, New York, 1974. (Reprinted by Dover, 2001).

\bibitem{ref3}
J. Arias de Reyna and J. van de Lune, \textit{A first encounter with the Riemann Hypothesis and its numerical verification}, La Gaceta de la RSME, \textbf{0} (1998), 1--24.

\bibitem{Conrey}
J. B. Conrey, \textit{More than two fifths of the zeros of the Riemann zeta function are on the critical line}, J. Reine Angew. Math., \textbf{399} (1989), 1--26.

\bibitem{levins}
N. Levinson, \textit{More than one third of zeros of Riemann's zeta-function are on $\sigma = 1/2$}, Advances in Mathematics, \textbf{13} (1974), 383--436.

\bibitem{ref32}
N. Levinson, \textit{Remarks on a formula of Riemann for his Zeta-function}, Journal of Mathematical Analysis and Applications, \textbf{41} (1973), 345--351.

\bibitem{ref19}
P. Borwein, S. Choi, B. Rooney, and A. Weirathmueller, \textit{The Riemann Hypothesis}, Springer, 2008.

\bibitem{ref59}
H. Rademacher, \textit{Topics in Analytic Number Theory}, Springer-Verlag, Berlin, 1973.

\bibitem{Selberg1}
A. Selberg, \textit{On the zeros of Riemann's zeta-function}, Skr. Norske Vid.-Akad. Oslo, No. 10, 1942.

\end{thebibliography}
\end{document}